# Optimal Battery Energy Storage Placement for Transient Voltage Stability Enhancement

Yongli Zhu[1], Chengxi Liu[2], Renchang Dai[1], Guangyi Liu[1], Yiran Xu[3]

[1]GEIRI North America, San Jose, USA (yongli.zhu@geirina.net, renchang.dai@geirina.net, guangyi.liu@geirina.net)
[2]Aalborg University, Aalborg, Denmark (cli@et.aau.dk)
[3]State Grid Nanjing Power Supply Company, Nanjing, Jiangsu, China

*Abstract*—**A placement problem for multiple Battery Energy Storage System (BESS) units is formulated towards power system transient voltage stability enhancement in this paper. The problem is solved by the Cross-Entropy (CE) optimization method. A simulation-based approach is adopted to incorporate higher-order dynamics and nonlinearities of generators and loads. The objective is to maximize the voltage stability index, which is set up based on certain grid-codes. Formulations of the optimization problem are then discussed. Finally, the proposed approach is implemented in MATLAB/DIgSILENT and tested on the New England 39-Bus system. Results indicate that installing BESS units at the optimized location can alleviate transient voltage instability issue compared with the original system with no BESS. The CE placement algorithm is also compared with the classic PSO (Particle Swarm Optimization) method, and its superiority is demonstrated in terms of fewer iterations for convergence with better solution qualities.**

*Index Terms*—**battery, cross entropy, energy storage, optimal placement, voltage stability**

## I. INTRODUCTION

Transient voltage stability, or short-term voltage stability [1], is the voltage issue involving dynamics of fast responsive loads such as induction motors, HVDC converters, etc. The study period of interest is in the order of several seconds and the analysis requires solution of system differential-algebraic-equations. This issue is usually tackled by installing SVC (Static Var Compensator) or STATCOM (Static Synchronous Compensator) [2]. In [3], a sensitivity index based on the system trajectory is utilized to determine the best location of SVC for voltage stability improvement. *VSI* (Voltage Stability Index) [4] based approaches are also proposed for choosing the best locations of var sources to alleviate the voltage instability issues [5], [6].

Besides SVC and STATCOM, another promising option to enhance the power system transient voltage stability is using BESS (Battery Energy Storage System). Nowadays, more and more BESS units have been deployed in power systems to provide services like peak shaving [7], system balancing [8], oscillation damping [9], [10], etc. In fact, the function of BESS can be further explored, i.e. not only to balance the system, but also to alleviate transient voltage issues.

One of the challenges to deploy the BESS for better system voltage stability is how to determine their locations in a power grid. Recent studies for optimal siting and sizing of BESS units in distribution network and microgrids [11]-[13] have been reported as well. Although most studies focused on the economic aspects for BESS placement problem, some efforts have been made for voltage stability improvement. In [14], a VIPI (Voltage Instability Proximity Index) based placement strategy is adopted for NaS battery in a distribution system from the perspective of static voltage stability, i.e. to ensure a large load margin. In [15], a TEF (Transient Energy Function) based optimal energy storage placement was studied in a microgrid. In [16], GA (Genetic Algorithm) is applied on the IEEE 14-bus system to search the best locations of SMES (Superconducting Magnetic Energy Storage) for voltage stability improvement using the *L*-index. Nevertheless, the placement of BESS for enhancing the transient voltage stability of bulk power system has not been fully investigated.

In this paper, an optimal BESS placement problem is formulated for transient voltage stability enhancement. It is then solved by the Cross-Entropy (CE) optimization method. The rest of the paper is organized as follows: Section II describes the power output model of BESS that is applied in this study; Section III presents the *VSI*-based problem formulation for BESS placement. Section IV introduces the principles of the CE optimization method. Section V did a case study on the New England 39-bus system. Section VI concludes this paper.

## II. MODELING OF BESS

### A. Overall Structrue of BESS

A BESS includes the battery cells and a converter interface called PCS (Power Conditioning System) as shown in Fig. 1. The BESS PCS typically includes a DC/DC converter for battery charging/discharging control and a DC/AC converter for interfacing the AC-grid. A battery cell can be modeled by an equivalent voltage source nonlinearly

This work is funded by State Grid Corporation Project SG5455HJ180021

depending on its SOC (State-Of-Charge), which is defined in (1), where $E_{total}$ is the energy capacity of BESS; [$t_0$, $t_1$] is a specific time interval for charging/discharging.

$$SOC = 1 - \Delta E / E_{total}, \quad \Delta E = \int_{t_0}^{t_1} P_{es}(t)dt \quad (1)$$

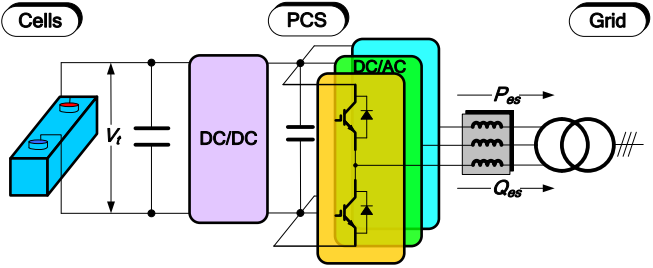

Figure 1. A typical topology of utility-scale BESS.

Most grid-related control strategies are implemented in the DC/AC part, and the cell voltage is maintained by the DC/DC part. Therefore, the BESS can be modeled as its PCS part while simultaneously considering the nonlinearity of the cell dynamics via its SOC limitations. In practice, the battery cell must be protected from deep charging or discharging for life-span consideration. In this study, the SOC range is assumed to be [20%, 80%].

### B. Power Output Model of BESS

For transient voltage stability studies, the BESS can be modeled by a first order transfer function shown in Fig. 2. The rationality of this modeling is that, the switching speed of the power electronic device is generally much faster than the electromechanical response of the synchronous generator. The effectiveness of such modeling approach for system level studies have been verified by previous researchers [17], [18].

Here, the BESS PCS regulates its active power output (positive when battery discharge) based on frequency deviation signals while keeping the PCS reactive power output at zero as show in Fig. 2. The terminal bus frequency of the BESS installation location is chosen as the input for $P_{ref}$.

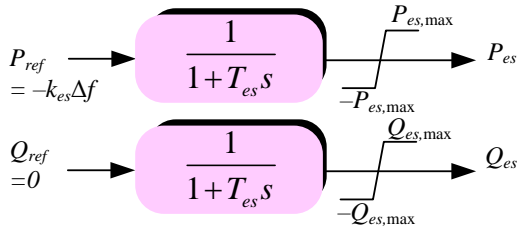

Figure 2. P-Q decoupled control scheme of the BESS PCS.

In Fig. 2, $k_{es}$ (positive) is the controller gain for each BESS, which stands for the theoretical maximum active power output of the BESS during transients and is set to 10.0 for each BESS in this paper. $T_{es}$ is the time constant of the BESS power converter, which is generally much smaller than the inertia time constants of large generators [18]. It is set to 0.02 sec in simulations. The energy capacity $E_{total}$ of a utility-scale BESS unit can range from 1 to 100 MWh [19]. In this paper, it is set to 10MWh. Finally, the BESS power output model is shown in (2) and (3).

$$P_{es} = \begin{cases} \dfrac{P_{ref}}{1+T_{es}s}, & SOC < SOC_{max} \text{ and } P_{ref} < 0 \\ 0, & SOC > SOC_{min} \text{ and } P_{ref} > 0 \\ 0, & \text{otherwise} \end{cases} \quad (2)$$

$$P_{es} \in [-P_{es,max}, P_{es,max}], \quad P_{ref} = -k_{es}\Delta f \quad (3)$$

## III. PLACEMENT PROBLEM FORMULATION

### A. Transient Voltage Stability Criteria

The NERC standards [20] define the post-fault voltage violations as illustrated in Fig. 3:

1) The post-fault instantaneous voltage dip or overshoot should be less than 25% for load buses and 30% for generator buses; any overshoot on or above 20% should be less than 20 cycles at load buses.

2) Post-transient voltage deviation should be less than 5% at all buses.

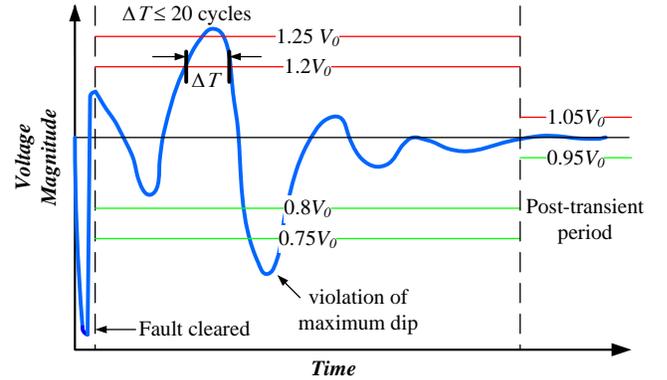

Figure 3. Post-fault transient voltage stability criteria for load bus.

### B. Voltage Stability Index

The Voltage Stability Index (*VSI*) used in this paper for the transient voltage stability study is based on [6]. Firstly, suppose there are top $k$ severe contingencies to be addressed. Then for bus $j$, define the voltage percentage variation at time step $t$ as ($N$ is the total number of buses and $j$ is the bus index):

$$D_{kj}^t = \left| \dfrac{V_{kj}^t - V_j^0}{V_j^0} \right|, j=1,\ldots,N, t=1,\ldots,T \quad (4)$$

Then, the following severity index (*SI*) can be defined:

$$SI_k = \dfrac{1}{T \times N} \sum_{t=1}^{T} \sum_{j=1}^{N} SI_{kj}^t, \quad SI_{kj}^t = \begin{cases} D_{kj}^t, & \textbf{criteria violated} \\ 0, & \textbf{otherwise} \end{cases} \quad (5)$$

The larger the $SI_k$, the severer the $k^{th}$ contingency can be. Then the next step is to setup a measure to quantify the contribution of BESS in terms of post-fault voltage improvement. In (6), a novel *maximum voltage recovery sensitivity* is used for *VSI* calculation. Note that, (6) is a generalized version of the *VSI* formula in [6], which defines the sensitivity for placing merely one var source each time.

$$VSI_{ij}^k = \frac{\max_t \{V_j^{k,new,t} - V_j^{k,old,t}, t=1,...,T\}}{\sum_{i=1}^{N_{es}} k_{es,i}} \quad (6)$$

In (6), the sum of $k_{es,i}$ for all the BESS power converters are used as denominator; thus, an interpretation of (6) is: the sensitivity *of voltage-recovery-extent* to the *theoretical maximum power outputs* of all placed BESS units.

The overall *VSI* for the $i^{th}$ bus in the $k^{th}$ contingency is:

$$VSI_i^k = \frac{1}{N}\sum_{j=1}^{N} VSI_{ij} \quad (7)$$

The overall *VSI* for the $i^{th}$ placement combination considering all the top-$K$ contingencies is:

$$VSI_i = \sum_{k=1}^{K} SI_k \cdot VSI_i^k \quad (8)$$

*C. Formulation*

Assume there are totally $N_{es}$ BESS units to be placed in a power grid to maximize the *VSI*, the optimal placement problem formulation for multiple BESS units is:

$$\max_{z_i} Obj = VSI(z) \quad (9)$$

$$\text{s.t.} \sum_{i=1}^{N} z_i = N_{es} \quad z_i \in \{0,1\} \quad (10)$$

In (10), $z_i$ is a binary variable with value "1" meaning that a BESS will be placed at the $i^{th}$ bus, and vice versa.

IV. CROSS ENTROPY OPTIMIZATION

The CE method was introduced by Rubinstein, *et. al.* [21]. Originally, the CE methods was developed for likelihood estimation of rare-events, i.e. events with very small probabilities (e.g. $\leq 10^{-4}$) It can be viewed as an "adaptive importance sampling" procedure that uses the cross-entropy (i.e. *Kullback-Leibler* (KL) divergence) as a measure of *closeness* between two sampling distributions. The principle of the method is iterated as follows:

1) Draw a sample based on a probabilistic distribution;

2) Based on the drawn sample, update parameters of that probabilistic distribution.

Thus, in each iteration, the candidate solution set is varying due to a *variable* probability density function.

*A. Basic Idea of CE Optimization*

Similar to the CE approach for rare event sampling, the original estimation-based CE method can be applied indirectly to the optimization problem. It assumes that the optimal solution is subjected to a *latent* probability distribution, which can be modeled or approximated by a function *f*.

Then, the seeking of optimal solution is mimicked by a consecutively sampling process of the "rare event"– **the event when the theoretical global optimal solution is found**. Specifically, for the following maximization problem (for minimization, the principle is similar), where $S$ is the objective function, **X** is the decision variables tread as random variable.

$$\gamma^* = \max S(\mathbf{X}) \quad (11)$$

under a density function $f(\cdot;u)$ (parametrized by $u$), the probability $l$ that $S(\mathbf{X})$ larger than some $\gamma$ is:

$$l = P(S(\mathbf{X}) \geq \gamma) = E_u[I_{(S(X)\geq \gamma)}] \quad (12)$$

where, $I_{(S(X)\geq\gamma)}$ (in short, "$I$" or "$I(\mathbf{X})$") is the "indicator function", with value "1" if the event in the bracket happens.

By using "Importance Sampling" technique [21], suppose $g$ is another probability function, and $g(x) = 0$ when $I = 0$. Then,

$$l = \int \left[I(x)\frac{f(x;u)}{g(x)}\right]g(x)\,dx = E_g[I(x)\frac{f(x;u)}{g(x)}] \Rightarrow$$
$$\hat{l} = \frac{1}{N}\sum_{i=1}^{N} I(x)\frac{f(\mathbf{X}_i;u)}{g(\mathbf{X}_i)}, g^* = \frac{f(x;u)I(x)}{\hat{l}} \quad (13)$$

Where, $g^*$ is the pdf under which the estimation variance of $l$ is minimal. The next step is to numerically find a $g$ as close to $g^*$ as possible. We can select $g = f(x; v)$ to approximate $f(x; u)$, with $v$ left to be solved. To this end, the *cross-entropy minimization* based on the "KL divergence" is utilized:

**KL divergence** $:= D(g^*, g) = E_{g^*}\left[\ln\frac{g^*(x)}{f(x)}\right], \mathbf{X} \sim g^*$

$$v^* = \arg\min_v D(g^*, f(x;v)) = \arg\max_v E_u I(\mathbf{X})\ln f(\mathbf{X};v) \quad (14)$$

where, $v$ is the new reference parameter used for the second-round importance sampling. The explicit solution of $v$ can be usually obtained when $f(x; v)$ is an *exponential-family* distribution like Gaussian distribution. Practically, for optimization problem, a multi-level CE method is adopted, i.e. choose $v^{(t)} = v^{(t-1)}$ in each iteration and $\gamma_t$ be the $(1-\rho)$ quantile of the sequence $S_1(X), ..., S_N(X)$ ($N$ is the sampling size), i.e.

$$P_{v^{(t-1)}}(S(X)\geq \gamma_t) \geq \rho \Leftrightarrow P_{v^{(t-1)}}(S(X)\leq \gamma_t)\leq 1-\rho \quad (15)$$

In (15), $\rho$ is a predefined ratio called "elite keeping rate", e.g. 50%. This process is iterated by updating $v^{(t-1)}$ to $v^{(t)}$ using (16) on $N_e=[\rho N]$ random samples that satisfy $S(\mathbf{X}) \geq \gamma_t$, until $v^{(t)}$ converges or maximum iterations reached.

$$v^{(t)} = \arg\max_v \frac{1}{N}\sum_{k=1}^{N} I(\mathbf{X}_k)\ln f(\mathbf{X}_k;v^{(t-1)}) \quad (16)$$

*B. CE method for Combinatorial Optimization*

For combinatorial optimization problem, the procedure of applying CE method is as follows:

1) Let the density function $f(\mathbf{x}; \mathbf{p})$ ($\mathbf{x}\in R^M$) be the *multivariate-Bernoulli* distribution with parameter $\mathbf{p}\in R^M$, i.e.

$$f(\boldsymbol{x};\boldsymbol{p}) = \prod_{j=1}^{M} p^{x_j}(1-p_j)^{(1-x_j)} \quad (17)$$

2) By (16), for $N$ samples $X_i \in \{0, 1\}^M$, $i = 1\ldots N$:

$$\hat{p} = \arg\max_{p} \frac{1}{N}\sum_{i=1}^{N} I_{\{S(X_i) \geq \hat{\gamma}_t\}} \ln f(X_i;p) \Rightarrow$$

$$\nabla \sum_{i=1}^{N} I_{\{S(X_i) \geq \hat{\gamma}_t\}} \ln f(X_i;p) = 0 \Leftrightarrow \hat{p}_j = \sum_{i=1}^{N} I_i x_{i,j} \bigg/ \sum_{i=1}^{N} I_i \quad (18)$$

Finally, the CE algorithm for combination optimization is summarized in Table I ($\alpha$ is a smoothing constant)

TABLE I.  CE ALGORITHM FOR COMBINATORIAL OPTIMIZATION

**Cross Entropy Algorithm for Combinatorial Optimization**

0  Input: $\rho$, $N$, $M$, $\alpha$, $N_e = \lceil \rho N \rceil$
1  Choose initial vector $\boldsymbol{p}^{(0)}$, $t=0$
2  Generate $N \times M$ sample matrix from **Ber** ($\boldsymbol{p}^{(t)}$) as:
   $X^{(t)} = [X_1^{(t)}, X_2^{(t)}, \ldots, X_M^{(t)}] \in R^{N \times M}$;
   $X_j^{(t)} = [x_{1,j}^{(t)}, x_{2,j}^{(t)}, \ldots, x_{N,j}^{(t)}]^T \in R^{N \times 1}$, $j=1,2,\ldots,M$
3  Evaluate and sort $S_i = S(X_i)$: $S_1 \leq \ldots \leq S_N$; updated $\gamma_t = S_{(N-N_e+1)}$; rearrange $X^{(t)}$ as $\tilde{X}^{(t)}$; Select the $N_e$ best samples
4  Update: for $j=1, 2\ldots, M$, $\hat{p}_j^{(t)} = \frac{1}{N_e}\sum_{i=1}^{N_e} \tilde{x}_{i,j}^{(t)}$
5  Smooth: $p_j^{(t)} = \alpha \hat{p}_j^{(t)} + (1-\alpha)\hat{p}_j^{(t-1)}$
6  If stop criteria reached, stop; otherwise, $t = t+1$, go to step-2

## V. CASE STUDY

In this section, the New England 10-machine, 39-bus system is chosen as the testbed for optimal BESS placement by the CE method.

### A. Network Modeling

The complete power system model is built in DIgSILENT/PowerFactory. The original 39-bus model in DIgSILENT is strong enough; thus, a weaker version is created by modifying the exciter parameters of some generators. The placed number of BESS is set to 3 in this study. The constant $k_{es}$ in (6) is set to 10.0 for each BESS.

### B. Simulation-based Optimizaiton Framework

The optimization algorithm is implemented in MATLAB using DIgSILENT as the simulation engine to dynamically evaluate objective values by the function call from MATLAB. The fault is applied at 0 sec and cleared after 0.1 sec at bus-16. Every time when a potential combination of BESS locations is generated during the CE method, it will be sent to DIgSILENT for simulation. Then the voltage waveforms will be returned to MATLAB for *VSI* calculation. In such simulation-based optimization framework, most nonlinearities of the power system models can be considered, thus the placement result will be more convincible.

### C. Placement Results

The bus set for the five severest contingencies ranked by *SI* values are {16,15,25,23,18}. The placement result is shown in Table II and Fig. 4. The selected locations and *VSI* values are presented. The iteration numbers needed for convergence is only 4. The CE parameters are: $\rho = 0.5$, $\alpha = 0.7$, $N=20$.

TABLE II.  PLACEMENT RESULT BY CE METHOD

| Siting Locations | VSI | Iter. No. at convergence | Max Iter. No. | Time Cost (s) |
|---|---|---|---|---|
| 34  35  36 | 20.7279 | 4 | 10 | 2522 |

### D. Comparison with PSO

The proposed BESS placement problem is also solved by the PSO (Particle Swarm Optimization) method with the *population size* set to 30. The result is listed in Table III. The iteration curves for both CE and PSO are shown in Fig. 5. PSO consumes more time but reaches a less-optimal objective value than the CE method.

TABLE III.  PLACEMENT RESULT BY PSO METHOD

| Siting Locations | VSI | Iter. No. at convergence | Max Iter. No. | Time Cost (s) |
|---|---|---|---|---|
| 23  34  35 | 20.6991 | 18 | 20 | 7057 |

From the results, it is obvious that, in this problem the CE method outperforms PSO in terms of less iteration numbers needed for convergence with similar objective values.

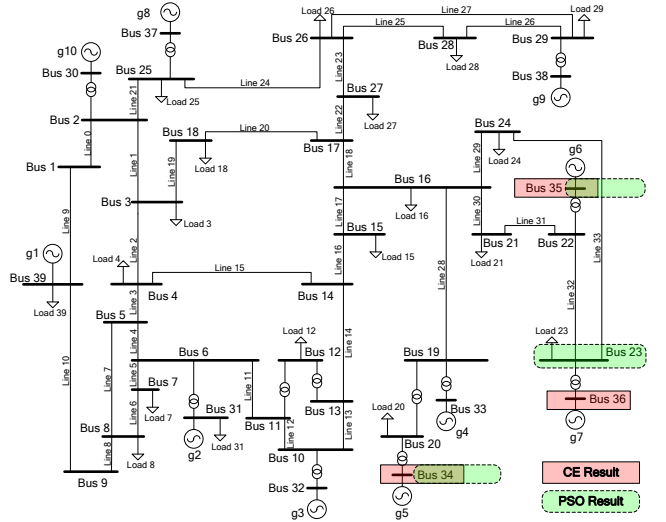

Figure 4.  BESS units placement on the New England 39-Bus system.

### E. Simulation Verifications

The optimal solution obtained is validated by time domain simulation for 5 seconds under the severest contingency of fault at bus-16. As an example, the voltage magnitudes of all the buses at 1.4 sec are displayed in Fig. 6. The overshot-event statistics are summarized in Table IV. It can be found that: without BESS, three buses (Bus 28, 29 and 38) exceed the

1.25 p.u. criteria; with BESS placed by PSO, only two buses violated; with BESS placed by CE, only one bus violated. In addition, the number of other buses who exceed the second criteria (i.e. 1.2 p.u.) is 2 for no BESS case; 1 for PSO case; and 0 for CE case. This validates the superiority of the placement result provided by CE method.

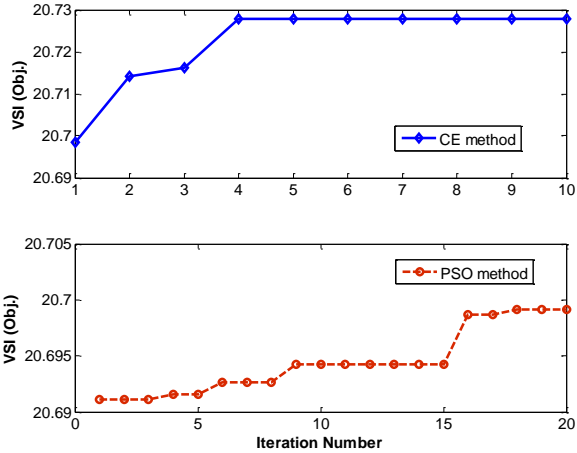

Figure 5. The iteration curves for CE and PSO methods respectively.

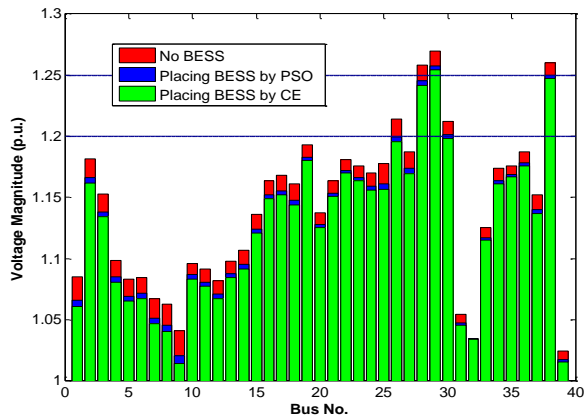

Figure 6. The post-fault voltage profile of the 39-bus system at 1.4 sec.

TABLE IV. STATISTICS OF VOLTAGE OVERSHOOT EVENTS AT TIME 1.4 SEC

|  | w/o BESS | BESS (PSO) | BESS (CE) |
|---|---|---|---|
| V >1.25 | 3 | 2 | 1 |
| 1.25 ≥V >1.2 | 2 | 1 | 0 |

## VI. CONCLUSION

The CE method has been successfully applied into the BESS placement problem. The placed BESS can enhance the transient voltage stability for the tested system. In our case, the CE method outperforms the PSO method in terms of an improved objective value and less iterations needed for convergence. Future work includes integrating the reactive-power control loop in the BESS and studying the effect of different parameters (sampling size, elite keeping rate, etc.) on the performance of the CE method.